\documentclass[dvips,twoside]{article}
\usepackage[latin1]{inputenc}
\usepackage[T1]{fontenc}
\usepackage{parskip}
\usepackage{amsmath,amsfonts,amssymb,amsxtra}
\usepackage{latexsym}
\usepackage{theorem}
\usepackage{graphicx,psfrag}
\usepackage{verbatim}

{\theoremstyle{plain}                       
\theorembodyfont{\itshape}
\newtheorem{Theorem}{Theorem}}

{\theoremstyle{plain}                       
\theorembodyfont{\itshape}
\newtheorem{Corollary}{Corollary}}

{\theoremstyle{plain}                       
\theorembodyfont{\itshape}
}

{\theoremstyle{plain}                       
\theorembodyfont{\rmfamily}
\newtheorem{Definition}{Definition}}

{\theoremstyle{plain}                       
\theorembodyfont{\rmfamily}
}

{\theoremstyle{plain}                       
\theorembodyfont{\rmfamily}
}

{\theoremstyle{plain}                       
\theorembodyfont{\rmfamily}
\newtheorem{Remark}{Remark}}

{\theoremstyle{plain}                       
\theorembodyfont{\itshape}
}

\begin{document}

\title{A necessary condition for the uniqueness of the stationary state of a Markov system}
\author{Ivan Werner\\
   {\small\it Email: ivan\_werner@pochta.ru}}
\maketitle

\begin{abstract}\noindent
 We continue the study of Markov systems started in \cite{Wer1}. In
 this paper, we prove a generalization of Breiman's strong low of
 large numbers \cite{Br} which implies a necessary condition for
 the uniqueness of the stationary state of a Markov system.

 \noindent{\it MSC}: 60J05, 37D35, 37A50, 37H99, 28A80

 \noindent{\it Keywords}:   Markov systems,
 iterated function systems with place-dependent probabilities, random systems with complete
connections, $g$-measures, Markov chains, equilibrium states,
fractals.
\end{abstract}

An introduction to Markov systems can be found in one of the
previous papers by the author, e.g. \cite{Wer1}, \cite{Wer3},
\cite{Wer6}, \cite{Wer7}, \cite{Wer8}.

Let $(K, d)$ be a complete separable metric space. Let
$\mathcal{M}:=(K_{i(e)}, w_e,p_e)_{e\in E}$ be a finite Markov
systems on $K$ (see figure 1) with an invariant Borel probability
measure $\mu$. We assume that $K_{i(e)}$ is open in $K=\bigcup_{e\in
E}K_{i(e)}$, $w_e|_{K_{i(e)}}$ is continuous and $p_e|_{K_{i(e)}}$
is continuous and bounded away from zero for all $e\in E$.
Obviously, these conditions imply the Feller property of the Markov
system. We shall denote by $C_B(K)$ the set of all bounded
continuous functions on $K$ and by $P(K)$ the set of all Borel
probability measures on $K$. Let $U:C_B(K)\longrightarrow C_B(K)$ be
the Markov operator associated with the Markov system, given by
\[Uf:=\sum\limits_{e\in E}p_ef\circ w_e\;\;\;\mbox{ for all }f\in
C_B(K),\] and $U^*:P(K)\longrightarrow P(K)$ be its adjoint
operator, given by
\[U^*\nu(f):=\int U(f)\ d\nu\;\;\;\mbox{ for all }f\in C_B(K),\
\nu\in P(K).\]

\begin{center}
\unitlength 1mm
\begin{picture}(70,70)\thicklines
\put(35,50){\circle{20}} \put(10,20){\framebox(15,15)}
\put(40,20){\line(2,3){10}} \put(40,20){\line(4,0){20}}
\put(50,35){\line(2,-3){10}} \put(5,15){$K_1$} \put(34,60){$K_2$}
\put(61,15){$K_3$} \put(31,50){\framebox(7.5,5)}
\put(33,45){\framebox(6.25,9.37)} \put(50,28){\circle{7.5}}
\put(45,21){\framebox(6,5)} \put(10,32.5){\line(6,1){15}}
\put(10,32.5){\line(3,-5){7.5}} \put(17.5,20){\line(1,2){7.5}}
\put(52,20){\line(2,3){4}} \put(13,44){$w_{e_1}$}
\put(35,38){$w_{e_2}$} \put(49,42){$w_{e_3}$}
\put(33,30.5){$w_{e_4}$} \put(30,15){$w_{e_5}$}
\put(65,37){$w_{e_6}$} \put(15,5){ Fig. 1. A Markov system.}
\put(0,60){$N=3$} \thinlines \linethickness{0.1mm}
\bezier{300}(17,37)(20,46)(32,52)
\bezier{50}(32,52)(30.5,51.7)(30,49.5)
\bezier{50}(32,52)(30,51)(28.7,51.7)
\bezier{300}(26,31)(35,36)(35,47)
\bezier{50}(35,47)(35,44.5)(33.5,44)
\bezier{50}(35,47)(35,44)(36,44) \bezier{300}(43,50)(49,42)(51,30)
\bezier{50}(51,30)(50.5,32)(49.2,32.6)
\bezier{50}(51,30)(50.6,32)(51.5,33.2)
\bezier{300}(39,20)(26,17)(18,25)
\bezier{50}(18,25)(19.5,24)(20,21.55)
\bezier{50}(18,25)(20,23.5)(22,24)
\bezier{300}(26,26)(37,28)(47,24) \bezier{50}(47,24)(45,25)(43,24)
\bezier{50}(47,24)(45,25)(44,26.5)
\bezier{100}(54.5,31.9)(56,37.3)(61,36.9)
\bezier{100}(61,36.9)(64.5,36.5)(66,34)
\bezier{100}(66,34)(68,30.5)(64.9,26.8)
\bezier{100}(64.9,26.8)(61.6,23.3)(57,23)
\bezier{50}(57,23)(58.5,23.3)(60.1,22.7)
\bezier{50}(57,23)(58.8,23.3)(59.5,24.8)
\end{picture}
\end{center}

We call a Markov system $\mathcal{M}$ {\it contractive} if and only
if there exists $0<a<1$ such that
\begin{equation*}
    \sum_{e\in E}p_{e}(x) d(w_{e}(x),w_{e}(y))\leq ad(x,y)\mbox{
for all } x,y\in K_{i(e)},\ e\in E.
\end{equation*}

 Let $\Sigma^+:=\{(\sigma_1,\sigma_2,...):\
\sigma_i\in E,\ i\in\mathbb{N}\}$ endowed with the product topology
of discreet topologies. For $x\in K$, let $P_x$ be the Borel
probability measure on $\Sigma^+$ given by
\[P_x( _1[e_1,...,e_n]):=p_{e_1}(x)p_{e_2}\circ
w_{e_1}(x)...p_{e_n}\circ w_{e_{n-1}}\circ ...\circ w_{e_1}(x),\]
for every cylinder set $_1[e_1,...,e_n]:=\{\sigma\in\Sigma^+:\
\sigma_1=e_1,...,\sigma_n=e_n\}$, which represents the Markov
process generated by the Markov system with the Dirac initial
distribution $\delta_x$.

It was shown in \cite{Wer7} that an irreducible contractive Markov
system $\mathcal{M}$ with uniformly continuous probabilities
$p_e|_{K_{i(e)}}$ has a unique invariant Borel probability measure
if $P_x\ll P_y$ for all $x, y\in K_{i(e)}$, $e\in E$, (this was
shown in \cite{Wer7} for some locally compact spaces, but it holds
also on complete separable spaces, as contractive $\mathcal{M}$ also
posses invariant measures on such spaces \cite{HS}).

This paper is motivated by the following question. Suppose the
Markov system $\mathcal{M}$ has a unique invariant Borel probability
measure. Does this imply some restrictions on the measures $P_x$,
$x\in K$? An answer to this question gives a generalization of
Breiman's strong low of large numbers \cite{Br}.

\begin{Theorem}\label{Bl}
 Suppose $1/n\sum_{k=0}^{n-1}U^kg(x)\to\mu(g)$
 for all $x\in K$, $g\in C_B(K)$. Let $f_e:K\longrightarrow [-\infty,+\infty]$ be
 Borel measurable such that $f_e|_{K_{i(e)}}$ is bounded and
 continuous for all $e\in E$. Then, for every $x\in K$,
 \[\lim\limits_{n\to\infty}\frac{1}{n}\sum\limits_{k=0}^{n-1}f_{\sigma_{k+1}}\circ
 w_{\sigma_k}\circ ...\circ w_{\sigma_1}(x)=\sum\limits_{e\in
 E}\int\limits_{K_{i(e)}}
 p_ef_e\ d\mu\] for $P_x$-a.e. $\sigma\in\Sigma^+$.
\end{Theorem}
{\it Proof.} We climb on the shoulders of Breiman \cite{Br}. By
Kakutani-Yosida norms ergodic lemma (e.g. see \cite{L} p. 441),
$1/n\sum_{k=0}^{n-1}U^kg$ converges uniformly to $\int g\ d\mu$ for
all $g\in C_B(K)$. Therefore, $1/n\sum_{k=0}^{n-1} U^k(p_ef_e)$
converges uniformly to $\int p_ef_e\ d\mu$ for all $e\in E$ (we use
the  convention $0\times\infty=0$).

Now, fix $x\in K$ and define
\begin{equation*}
    X_n(\sigma):=\left\{\begin{array}{cc}
    w_{\sigma_n}\circ w_{\sigma_{n-1}}\circ...\circ w_{\sigma_1}(x), & n\geq 1 \\
    x,& n=0,
  \end{array}\right.
\end{equation*}
\[X'_n(\sigma):=f_{\sigma_n}\circ X_{n-1}(\sigma)\]
for all $\sigma\in\Sigma^+$, $n\in\mathbb{N}$,
\begin{equation*}
    Z^1_n:=\left\{\begin{array}{cc}
    X'_n-E\left(X'_n|X_{n-1}\right), & n> 1 \\
    0,& n\leq 1,
  \end{array}\right.
\end{equation*}
and
\begin{equation*}
    Z^k_n:=\left\{\begin{array}{cc}
    E\left(X'_n|X_{n-k+1}\right)-E\left(X'_n|X_{n-k}\right), & n> k \\
    0,& n\leq k
  \end{array}\right.
\end{equation*} for $k>1$, where E(.|.) denotes the conditional
expectation with respect to measure $P_x$.

We are going to use the following result (\cite{L}, p. 387). Let
$Y_1,Y_2,...$ be a sequence of random variables such that
$E(Y_n|Y_{n-1},...,Y_1)=0$ and $EY_n^2\leq c<\infty$ for all $n$.
Then $1/m\sum_{l=1}^mY_l\to 0$ a.s..

To apply this, note that
\[E(Z^k_n|Z^k_{n-1},...,Z^k_1)=E(E(Z^k_n|X_{n-k},X_{n-k-1},...,X_1)|Z^k_{n-1},...,Z^k_1),\]
and that, since the $X_1,X_2,...$ form a Markov chain,
\[E(Z^k_n|X_{n-k},X_{n-k-1},...,X_1)=E(Z^k_n|X_{n-k})=0.\]
Furthermore, $E(Z^k_n)^2\leq 2\max_{e\in
E}||f_e|_{K_{i(e)}}||^2_{\infty}$. Thus,
\begin{equation}\label{sl}
    \frac{1}{m}\sum\limits_{l=1}^mZ^k_l\to 0\;\;\; P_x\mbox{-a.e.}.
\end{equation}
Now, write
\[X'_n-E(X'_n|X_{n-k})=Z^1_n+Z^2_n+...+Z^k_n,\;\;\; n>k.\]
Thus, by (\ref{sl}),
\[\lim\limits_{m\to\infty}\left|\frac{1}{m}\sum\limits_{n=1}^mX'_n-
\frac{1}{m}\sum\limits_{n=k+1}^mE(X'_n|X_{n-k})\right|=0\;\;\;
P_x\mbox{-a.e.}.\] Or, neglecting at most $k$ terms,
\[\lim\limits_{m\to\infty}\left|\frac{1}{m}\sum\limits_{n=1}^mX'_n-
\frac{1}{m}\sum\limits_{n=1}^mE(X'_{n+k}|X_{n})\right|=0\;\;\;
P_x\mbox{-a.e.}.\] So that, for fixed $s$,
\begin{equation}\label{pl}
    \lim\limits_{m\to\infty}\left|\frac{1}{m}\sum\limits_{n=1}^mX'_n-
\frac{1}{m}\sum\limits_{n=1}^m\left(\frac{1}{s}\sum\limits_{k=1}^sE(X'_{n+k}|X_{n})\right)\right|=0\;\;\;
P_x\mbox{-a.e.}.
\end{equation}
Now, observe that
\begin{eqnarray*}
&&E(X'_{n+k}|X_n)\\
&=&\sum\limits_{e_{n+1},...,e_{n+k}}p_{e_{n+1}}(X_n)p_{e_{n+2}}(X_n)...\\
\;\;\;\;\;\;&&\times p_{e_{n+k}}\circ w_{e_{n+k-1}}\circ...\circ
w_{e_{n+1}}(X_n)f_{e_{n+k}}\circ
w_{e_{n+k-1}}\circ...\circ w_{e_{n+1}}(X_n)\\
&=&\sum\limits_{e_{n+k}}U^{k-1}(p_{e_{n+k}}f_{e_{n+k}})(X_n).
\end{eqnarray*}
Since $1/s\sum_{k=1}^s U^{k-1}(p_ef_e)$ converges uniformly to $\int
p_ef_e\ d\mu$ for all $e\in E$, for every $\epsilon>0$, we can
choose $s>0$ such that
\[\|\frac{1}{s}\sum\limits_{k=1}^s\sum\limits_{e\in
E}U^{k-1}(p_ef_e)-\sum\limits_{e\in E}\int p_ef_e\
d\mu\|_\infty<\epsilon.\] By (\ref{pl}), for such $s$,
\[\lim\limits_{m\to\infty}\left|\frac{1}{m}\sum\limits_{k=1}^mX'_n-\sum\limits_{e\in
E}\int p_ef_e\ d\mu\right|\leq\epsilon\;\;\; P_x\mbox{-a.e.}.\]
\hfill$\Box$

\begin{Corollary}\label{fc}
 Suppose $K$ is a compact metric space and $\mu$ is a unique
 invariant Borel probability measure of the Markov system
 $\mathcal{M}$. Let $f_e: K\to [-\infty,+\infty]$ be Borel
 measurable such that $f_e|_{K_{i(e)}}$ is continuous for all $e\in
 E$. Then, for every $x\in K$,
 \[\lim\limits_{n\to\infty}\frac{1}{n}\sum\limits_{k=0}^{n-1}f_{\sigma_{k+1}}\circ
 w_{\sigma_k}\circ ...\circ w_{\sigma_1}(x)=\sum\limits_{e\in
 E}\int\limits_{K_{i(e)}}p_ef_e\ d\mu\] for $P_x$-a.e. $\sigma\in\Sigma^+$.
\end{Corollary}
{\it Proof.} By the hypothesis $P(K)$ is weakly$^*$ compact. Let
$x\in K$. Since every weakly$^*$ convergent subsequence
$1/n_m\sum_{k=0}^{n_m-1}{U^*}^k\delta_x$ converges to an invariant
Borel probability measure of $U^*$, by the uniqueness of $\mu$, we
conclude that
\[\frac{1}{n}\sum\limits_{k=0}^{n-1}{U^*}^k\delta_x\stackrel{w^*}{\to}\mu.\]
Hence, the claim follows by Theorem \ref{Bl}.\hfill$\Box$

\begin{Remark}
   Obviously, Corollary \ref{fc} is a generalization of Breiman's
   strong low of large numbers \cite{Br} for finite Markov systems.
\end{Remark}

The following is a version for contractive Markov systems.
\begin{Corollary}\label{sc}
 Suppose $K$ is a complete separable metric space. Suppose $\mathcal{M}$ is a contractive Markov system
  with a unique invariant Borel probability measure $\mu$. Let $f_e: K\to [-\infty,+\infty]$ be Borel
 measurable such that $f_e|_{K_{i(e)}}$ is bounded and continuous for all $e\in
 E$. Then, for every $x\in K$,
 \[\lim\limits_{n\to\infty}\frac{1}{n}\sum\limits_{k=0}^{n-1}f_{\sigma_{k+1}}\circ
 w_{\sigma_k}\circ ...\circ w_{\sigma_1}(x)=\sum\limits_{e\in
 E}\int\limits_{K_{i(e)}}p_ef_e\ d\mu\] for $P_x$-a.e. $\sigma\in\Sigma^+$.
\end{Corollary}
{\it Proof.} It was shown in \cite {HS} that the hypothesis implies
\[\frac{1}{n}\sum\limits_{k=0}^{n-1}U^kg(x)\to\mu(g)\]
for all $x\in K$, $g\in C_B(K)$. Thus, the claim follows by Theorem
\ref{Bl}.\hfill$\Box$

\begin{Definition}
    Let $\Sigma:=\{(...,\sigma_{-1},\sigma_0,\sigma_1,...):e_i\in E\
\forall i\in\mathbb{Z}\}$ endowed with the product topology of
discreet topologies. Denote by $S$ the left shift map on $\Sigma$.
We call the shift invariant Borel probability measure $M$ on
$\Sigma$ given by
\[M\left(_m[e_1,...,e_k]\right):=\int
p_{e_1}(x)p_{e_2}(w_{e_1}x)...p_{e_k}(w_{e_{k-1}}\circ...\circ
w_{e_1}x)d\mu(x),\] for every cylinder set
$_m[e_1,...,e_k]:=\{\sigma\in\Sigma:\ \sigma_i=e_i\mbox{ for all
}m\leq i\leq n\}$, the {\it generalized Markov measure} (see
\cite{Wer3}). We call the measure preserving transformation
$S:(\Sigma,M)\longrightarrow (\Sigma,M)$ a {\it generalized Markov
shift}.
\end{Definition}

The following corollary gives a necessary condition on measures
$P_x$, $x\in K$, for the uniqueness of $\mu$.
\begin{Corollary}\label{tc}
 Suppose $\mu$ is a unique invariant Borel probability measure of
 $\mathcal{M }$. Suppose $K$ is a compact metric space, or
 $K$ is a complete separable metric space and $\mathcal{M}$ is contractive.  Then, for every $x\in K$,
 \[\lim\limits_{n\to\infty}\frac{1}{n}\log P_x( _1[\sigma_1,...,\sigma_n])=\sum\limits_{e\in
 E}\int\limits_{K_{i(e)}} p_e\log p_e\ d\mu\]
 for $P_x$-a.e. $\sigma\in\Sigma^+$.
\end{Corollary}
{\it Proof.}  Set $f_e:=\log p_e$ for all $e\in E$. Then, by
Corollary  \ref{fc} or Corollary \ref{sc}, for every $x\in K$,
\[\lim\limits_{n\to\infty}\frac{1}{n}\log [p_{\sigma_1}(x)p_{\sigma_2}\circ w_{\sigma_1}(x)...
p_{\sigma_n}\circ w_{\sigma_{n-1}}\circ...\circ
w_{\sigma_1}(x)]=\sum\limits_{e\in E}\int\limits_{K_{i(e)}} p_e\log
p_e\ d\mu\] for $P_x$-a.e. $\sigma\in\Sigma^+$.\hfill$\Box$

\begin{Remark}
    Note that $-\sum_{e\in E}\int_{K_{i(e)}} p_e\log p_e\ d\mu$ is the
    Kolmogorov-Sinai entropy of the generalized Markov shift
    associated with $\mathcal{M}$ and $\mu$ (this was proved in \cite{Wer8} under additional
    assumption of uniform continuity of each $p_e|_{K_{i(e)}}$). This
    entropy formula also plays a central role in \cite{S}.
\end{Remark}

\end{document}